\documentclass[11pt,a4paper]{article}
\usepackage[T1]{fontenc}
\usepackage[utf8]{inputenc}
\usepackage{authblk}
\usepackage{algorithm}
\usepackage{algorithmic}
\usepackage{paralist}
\usepackage{amsmath,amssymb}
\usepackage{fancyhdr}
\usepackage{hyperref}
\usepackage{pgf}
\usepackage{tikz} \usetikzlibrary{automata, arrows}

\usepackage{pgfplots}
\pgfplotsset{compat=newest}
\usepgfplotslibrary{units}

\usetikzlibrary{arrows,positioning,calc}
\usepackage{tikz}
\usepackage{tkz-graph}
\usetikzlibrary{backgrounds}
\usepgflibrary{shapes}
\usetikzlibrary{through}
\usepgflibrary{shapes.symbols} 
\usepgflibrary[shapes.symbols] 
\usepgflibrary{shapes.symbols}
\usetikzlibrary{shapes.symbols} 
\usetikzlibrary[shapes.symbols] 
\usepgflibrary{shapes.geometric} 
\usepgflibrary[shapes.geometric] 
\usetikzlibrary{shapes.geometric} 
\usetikzlibrary[shapes.geometric] 
\usetikzlibrary{automata,positioning}
\usetikzlibrary{shapes,snakes}

\tikzset{
    >=stealth',
    smallpunkt/.style={
               rectangle,
               rounded corners,
               draw=black, very thick,
               text width=1.5em,
               minimum height=2em,
               text centered},
    punkt/.style={
           rectangle,
           rounded corners,
           draw=black, very thick,
           text width=6.5em,
           minimum height=2em,
           text centered},
     largepunkt/.style={
                rectangle,
                rounded corners,
                draw=black, very thick,
                text width=9.5em,
                minimum height=2em,
                text centered},
    pil/.style={
           ->,
           thick,
           }
}

\pgfdeclarelayer{background}
\pgfdeclarelayer{foreground}
\pgfsetlayers{background,main,foreground}

\usepackage{epsfig}
\usepackage{amsthm}

\def\squarebox#1{\hbox to #1{\hfill\vbox to #1{\vfill}}}

\newcommand{\dahntab}[1]{
  \newbox\mybok%
  \setbox\mybok=\hbox{\vbox{
      \begin{tabbing}
        #1
      \end{tabbing}%
    }}

  \newdimen\bokwidth%
  \bokwidth=\wd\mybok%
  \newdimen\myl%
  \myl=\textwidth%
  \divide\myl by 2%
  \divide\bokwidth by -2%
  \advance\myl by\bokwidth%
  \vrule width\myl height 0pt depth 0pt%
  \usebox\mybok%
}

\newcommand{\remove}[1]{}
\newcommand{\dx}[1]{\,\mathrm{d} #1}
\newcommand{\fd}[1]{\mathrm{D}_*^{\alpha}\, #1}
\newcommand{\ceil}[1]{\left\lceil #1 \right\rceil}
\newcommand{\floor}[1]{\left\lfloor #1 \right\rfloor}

\begin{document}

\title{Parallel simulations for Fractional-Order Systems}
\author[1]{Andrada BABAN}
\author[1]{Cosmin BONCHI\c{S}}
\author[2]{Alexandru FIKL}
\author[1]{Florin RO\c{S}U}
\affil[1]{West University of Timi\c{s}oara and the eAustria Research Institute, Bd. V. P\^arvan 4, cam 045B, Timi\c{s}oara, RO-300223, Romania.} 
\affil[2]{Dept. of Aerospace Engineering, University of Illinois  at Urbana-Champaign, USA.}

\def\titlerunning{Parallel simulations for Fractional-Order Systems}
\def\authorrunning{A. BABAN, C. BONCHI\c{S}, A. FIKL \& F. RO\C{S}U}

\maketitle

\begin{abstract}

In this paper, we explore how numerical calculations can be accelerated by implementing several numerical methods of fractional-order systems using parallel computing techniques. We investigate the feasibility of parallel computing algorithms and their efficiency in reducing the computational costs over a large time interval. Particularly, we present the case of Adams-Bashforth-Mouhlton predictor-corrector method and measure the speedup of two parallel approaches by using GPU and HPC cluster implementations.
\end{abstract}

{\bf Keywords:} Fractional-order systems, parallel numerical algorithms, GPU processing, HPC processing

\section{Introduction}
It is well understood that fractional-order derivatives provide a good tool for the description of memory and hereditary properties of various processes,  fractional-order systems being characterized by infinite memory. Generalizations of dynamical systems using fractional-order derivatives instead of classical integer-order derivatives have proved to be useful and more accurate in the mathematical modeling of real world phenomena arising from several interdisciplinary areas such as: diffusion and wave propagation, viscoelastic liquids, fractional kinetics, boundary layer effects in ducts, electromagnetic waves, electrode-electrolyte polarization.

Theoretical characterization of chaos in fractional-order dynamical systems is yet to be investigated. However, chaotic behavior has been observed by numerical simulations in many systems such as: a fractional-order Van der Pol system \cite{Barbosa2007}, fractional-order Chua and Chen's systems \cite{Donato2008Bif, Donato2008Frac}, a fractional-order Rossler system \cite{Weiwei2009} and a fractional-order financial system \cite{Wang2010}. Nevertheless, it is worth noting that numerical simulations are limited by the fact that they only reveal the chaotic behavior of discrete-time dynamical systems that are obtained by discretizing the fractional-order systems. 

In order to assess chaotic behavior of fractional-order dynamical systems, it is extremely important to be able to accurately estimate the solutions over a large time interval. Several numerical methods are used for fractional-order systems, such as a generalization of the Adams-Bashforth-Moulton \cite{Diethelm2002} predictor-corrector method  or a class \cite{Galeone2009} of p-fractional linear multistep methods. 

For the numerical assessment of chaotic behavior of fractional-order systems, it is extremely important to accurately estimate the solutions over a large time interval. In numerical schemes such as the Adams-Bashforth-Moulton predictor-corrector method, due to the hereditary nature of the problem and the non-locality of fractional-order derivatives, at every iteration step all previous iterations have to be taken into account. It may be feasible to overcome these difficulties with the aid of parallel computing algorithms (e.g. \cite{Diethelm2011}) implemented in a conventional way or using available high performance computing systems. 

In this paper we present the case of Adams-Bashforth-Moulton predictor-corrector method and measure the speedup of two parallel approaches by using GPU and HPC implementations. In the following, we will present the theoretical problem, we will describe the parallel dynamic scheme and then we will present our simulation results and our conclusions. 

\section{Preliminaries}
Were trying to solve an ordinary fractional differential equation of the form:
\begin{equation}\label{eq:fde}
\begin{cases}
\fd{y(t)} = f(t, y(t)), & t \in [0, T]\\
y^{(k)}(0) = y_0^{k}, & k \in \{0, \dots, \ceil{\alpha} - 1\},
\end{cases}
\end{equation}
where $\alpha \in [0, 1]$ and $\ceil{\cdot}$ denotes the ceiling function that
rounds up to the nearest integer. The fractional derivative is defined as:
\[
\fd{y(t)} = \frac{1}{\Gamma(\ceil{\alpha} - \alpha)}
            \int_0^T \frac{y^{(\ceil{\alpha})}(t)}{(t - \tau)^{\alpha - \ceil{\alpha}+ 1}} \dx{\tau}.
\]

The numerical method used to solve~\eqref{eq:fde} is a fractional version of the
Adams-Bashforth-Moulton method (a classic predictor-collector method) presented also in Figure \ref{alg:ABM}. The
domain $[0, T]$ is discretized into $N$ intervals with a step size $h = \frac{T}{N}$
and the grid points $t_n = nh$, for $n \in \{0, \dots, N\}$. We will also denote
$y_n = y(t_n)$ and $f_n = f(t_n, y_n)$ with $y_0$ as the initial condition.
Since we restrict $\alpha$ to the $[0, 1]$ interval, we do not require any of
the higher derivatives of the initial condition.

The first step of the scheme is the \textbf{predictor}, which will give a
first approximation $y_{n + 1}^P$ of our solution:
\begin{equation}\label{eq:abm}
y_{n + 1}^P = \sum_{k = 0}^{\ceil{\alpha} - 1} \frac{t_{n + 1}^k}{k!} y_0^{(k)} +
              h^\alpha \sum_{k = 0}^n b_{n - k} f_k,
\end{equation}
where the weights $b_{n}$ are given by:
\[
b_n = \frac{(n + 1)^\alpha + n^\alpha}{\Gamma(\alpha + 1)}.
\]

The second and final approximation of our solution, called the \textbf{corrector},
is given by:
\begin{align*}
y_{n + 1} & = \sum_{k = 0}^{\ceil{\alpha} - 1} \frac{t_{n + 1}^k}{k!} y_0^{(k)} + h^\alpha \left(c_n f_0 + \sum_{k = 1}^n a_{n - k} f_k + \frac{f(t_{n + 1}, y^P_{n + 1})} {\Gamma(\alpha + 2)}\right),
\end{align*}

where the weights $a_n$ and $c_n$ are defined as:
\[
a_n = \frac{(n + 2)^{\alpha + 1} - 2(n + 1)^{\alpha + 1} + n^{\alpha + 1}}{\Gamma(\alpha + 2)}
\]
and
\[
c_n = \frac{n^{\alpha + 1} - (n - \alpha)(n + 1)^\alpha}{\Gamma(\alpha + 2)}.
\]

The most difficult part of numerically solving such equations consists in the
fact that at each step, we require the complete history of the variable. That
is to say, to compute $y_{n + 1}$, we need to know all the previous values
$y_k$ that are used to compute $f_k$, for $k \leq n$. This makes fractional
differential equations notoriously hard to parallelize.

\section{Dynamic Scheme}

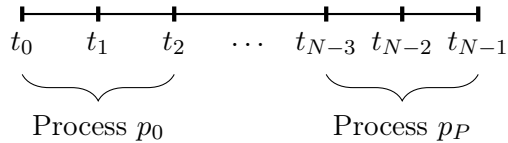
\begin{figure}[!ht]
\centering
\begin{tikzpicture}[scale=2]
\draw [thick] (0, 0) -- (3, 0);
\foreach[count=\i from 0] \x in {0, 0.5,...,1}{
    \draw[ultra thick] (\x, -0.05) -- (\x, 0.05);
    \node at (\x, -0.05) [below] {$t_{\i}$};
}
\foreach[count=\i from 0]\x in {2, 2.5,...,3}{
    \pgfmathtruncatemacro{\n}{3 - \i};
    \draw[ultra thick] (\x, -0.05) -- (\x, 0.05);
    \node at (\x, -0.05) [below] {$t_{N - \n}$};
}

\draw [decorate,decoration={brace,amplitude=10pt,mirror}] (0,-0.4) -- (1,-0.4);
\draw [decorate,decoration={brace,amplitude=10pt,mirror}] (2,-0.4) -- (3,-0.4);

\node at (0.5, -0.6) [below] {Process $p_0$};
\node at (2.5, -0.6) [below] {Process $p_P$};
\node at (1.5, -0.2) {$\dots$};
\end{tikzpicture}
\caption{Partition of time steps between processes.}
\label{fig:partition}
\end{figure}

The first algorithm we are going to look at proposes to break up all the $N$
time steps between the $P$ processes we have, giving: $ \overline{N} = \frac{N}{P}$
steps per process (see Figure~\ref{fig:partition})

During an iteration $n$, a process $p$ that represents the block in which
we can find $n$, that is to say if
\[
n \in \left[\overline{N} \cdot p, \overline{N} \cdot (p + 1)\right]
\]
will compute the value of $y_n$. This is done by computing a partial sum in
each process $p' < p$ and gathering the results to finally form the complete
sum as defined in~\eqref{eq:abm}.

One major downside of the algorithm we are going to describe is that given a
fix number $P$ of processes, most of them will be idle for a big part of the
computations. To be exact, they will be idle as long as:
\[
\floor{\frac{n}{\overline{N}}} < p,
\]
that is to say, as long as the current iteration has not reached the block
corresponding to process $p$. 

\begin{figure}[h]
\begin{center}
\begin{minipage}[c]{0.7\columnwidth}
\begin{center}
\dahntab{
\=\ \ \ \ \=\ \ \ \ \=\ \ \ \ \= 
{\bf ParallelABM Algorithm}
\\
\\
{\bf INPUT} : ParABM$(T, N, P, p, y_0)$ \\
\>	$T$ end of the time interval. \\
\>    $N$ global number of points. \\
\>    $P$ number of processes. \\
\>    $p$ current process. \\
\>    $y_0$ initial condition\\    
\>    $\overline{N} := N / P$ \\
\>    $n_{min} := \overline{N} \cdot p$\\
\>    $n_{max} := \overline{N} \cdot (p + 1)$\\
\>    {\bf for} $n \in [1, N]$\\
\> \>        $S := 0$\\
\> \>        $p' := \floor{\frac{n}{\overline{N}}}$ process computing $y_n$\\
\> \>        {\bf if} $p > p'$ {\bf then}\\
\> \> \>             continue\\
\> \>         {\bf end}\\
\> \>         {\bf if} $p = p'$ {\bf then}\\
\> \> \>            receive partial sums from previous processes\\
\> \> \>             {\bf for} $p' \in [0, p - 1]$ \\
\> \> \> \>                 MPI\_Recv$(p', S_{p'})$\\
\> \> \>   \>               $S := S + S_{p'}$\\
\> \> \>             {\bf end}\\
\> \> \>            compute local sum\\
\> \> \>         {\bf for} $k \in [n_{min}, n - 1]$ \\
\> \> \> \>                $S := S + b_{n - k} f_k$\\
\> \> \>            {\bf end}\\
\> \> \>            compute the predictor at time $t_n$\\
\> \> \>             $y_n^P := y_0 + S$\\
\> \>         {\bf else}\\
\> \> \>            compute the partial sum and send it to $p'$\\
\> \> \>            {\bf for} $k \in [n_{min}, n_{max}]$\\
\> \> \> \>                $S := S + b_{n - k} f_k$\\
\> \> \>            {\bf end}\\
\> \> \>            MPI\_Send$(p', S)$\\
\> \>        {\bf end}\\
\> \>         repeat the same steps to compute the corrector $y_{n}$\\
\>  {\bf end}\\
{\bf OUTPUT} : Imputation $Y=(y_{i})_{i\in [1,N]}$.
}
\end{center}
\end{minipage}
\end{center}
\caption{\label{alg:ABM}Parallel Algorithm for the Adams-Bashforth-Moulton (ABM) scheme.}
\end{figure}

\section{Simulation results}

\subsection*{BlueGene simulations}

 On our tests we used UVT’s BlueGene/P cluster that consists of an fully loaded single BlueGene/P rack that has move than 1000 CPUs and 4TB of RAM memory. It can offer a 11.7 TFlops sustained performance.
	One important issue of the algorithm described previously (Algorithm \ref{alg:ABM}) is that for a fix number $P$ of processes, most of them will be idle for the time when the sum is computed from the partial results that are gathered.

	Below are several simulation to find a balance for this drawback. We made several simulation of computing a number of elements, with different number of processes, and we observe the average running time for different time computation per process:
\vspace{-0.1cm}

\begin{table}[h]
\caption{BlueGene simulation} 
\centering 
\begin{tabular}{r r r r r}
\hline                        
 Number of  & Average& time~~~ per& processes (s)\\  [0.5ex] 
 steps (N)  & 64 & 128 & 256 \\  [0.5ex] 
\hline                 
100000 & 163.08 & 100.23 & 194.84\\
500000 & 1146.45 & 1889.18 & 1566.70 \\
1000000 & 660.61 & 1871.48 & 326.51 \\
 1500000  & - & 1628.62 & 572.29\\
 2000000  & - & 526.30 & 1829.73\\
\hline 
\end{tabular}
\label{table:BlueGene} 
\end{table}

\vspace{-0.3cm}
For a certain number of elements, for example for 1 million elements and 64 processes, the average time is 660 seconds of computing time. By double the number of processes, this leads to fewer elements to compute per process, but increase the time to compute the final result from the partial sums, and the time is almost triple. But increasing more, 4 times, it reduces the time by half. The explanation for this is, although on each iteration the time for computing the final result is increasing, the advantage is decisive that there are fewer elements per process to compute.

The final goal is to compute as much as possible elements in a reasonable time frame. And the times for 1.5 millions and 2 millions are very interesting. Depending on the number of processes, the computing times are swapped between 128 and 256 processes.
Further investigations will provide more data to identify a pattern and predict the right combination to minimize the computing time for a certain number of elements.

\vspace{-0.3cm}
\begin{table}[ht]
\caption{Average running time results} 
\centering 
\begin{tabular}{c r r r}
\hline                        
Test & Number of  & Average HPC & Average CUDA\\  [0.5ex] 
case & steps (N)  & running time & running time\\  [0.5ex] 
\hline                 
1 & 100000 & 170.75 & 169.06\\
2 & 500000 & 1592.66& 1276.42 \\
3 & 1000000 & 4621.25 & 2670.56\\
4 &  1500000  & 9162.33 & 4450.09\\
5 &  2000000 & 14931.16 & 6639.45 \\
6 & 2500000 & 22697.66 & 9229.66\\
7 & 3000000 & 31659.66 & 12226.39\\ [1ex]      
\hline 
\end{tabular}
\label{table:avgRunTime} 
\end{table}	

\subsection*{CUDA simulations}

The Parallel Algorithm for the Adams-Bashforth-Moulton scheme using MPI (Algorithm 1) was been also adapted for CUDA. The ParABM algorithm in this case suffered few changes: the partial sums are computed using parallel sum reduction; the weights $a_n$, $b_n$, $c_n$ not being dependent on terms being computed in advance via previous steps, they can be preprocessed in parallel, taking into consideration that each thread n computes the weight at time n.

\begin{figure} [h]
\begin{center}
\includegraphics[width=0.65\textwidth]{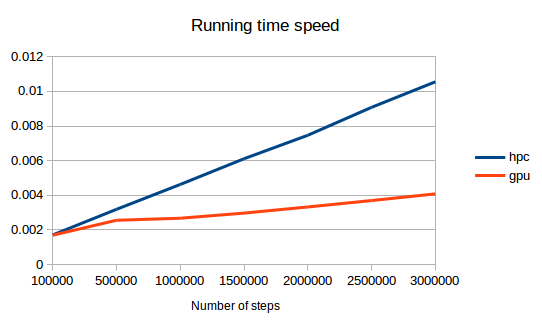}
\end{center}
\caption{Running time speed HPU CPU vs GPU}
\label{fig:hpcvsgpu}
\end{figure}

In Table \ref{table:avgRunTime} we compare all together the absolute average time of multiple GPU cores versus multiple CPUs. We can observe in Figure \ref{fig:hpcvsgpu} that in both cases the absolute processing time is increasing with number of steps. It is surprising that we obtained on HPC CPU a lower rate versus GPU, that tell us we should look carefully to the HPC processing time and idle time on the future. 

\begin{figure}[h] 
\begin{center}
\includegraphics[width=0.40\textwidth]{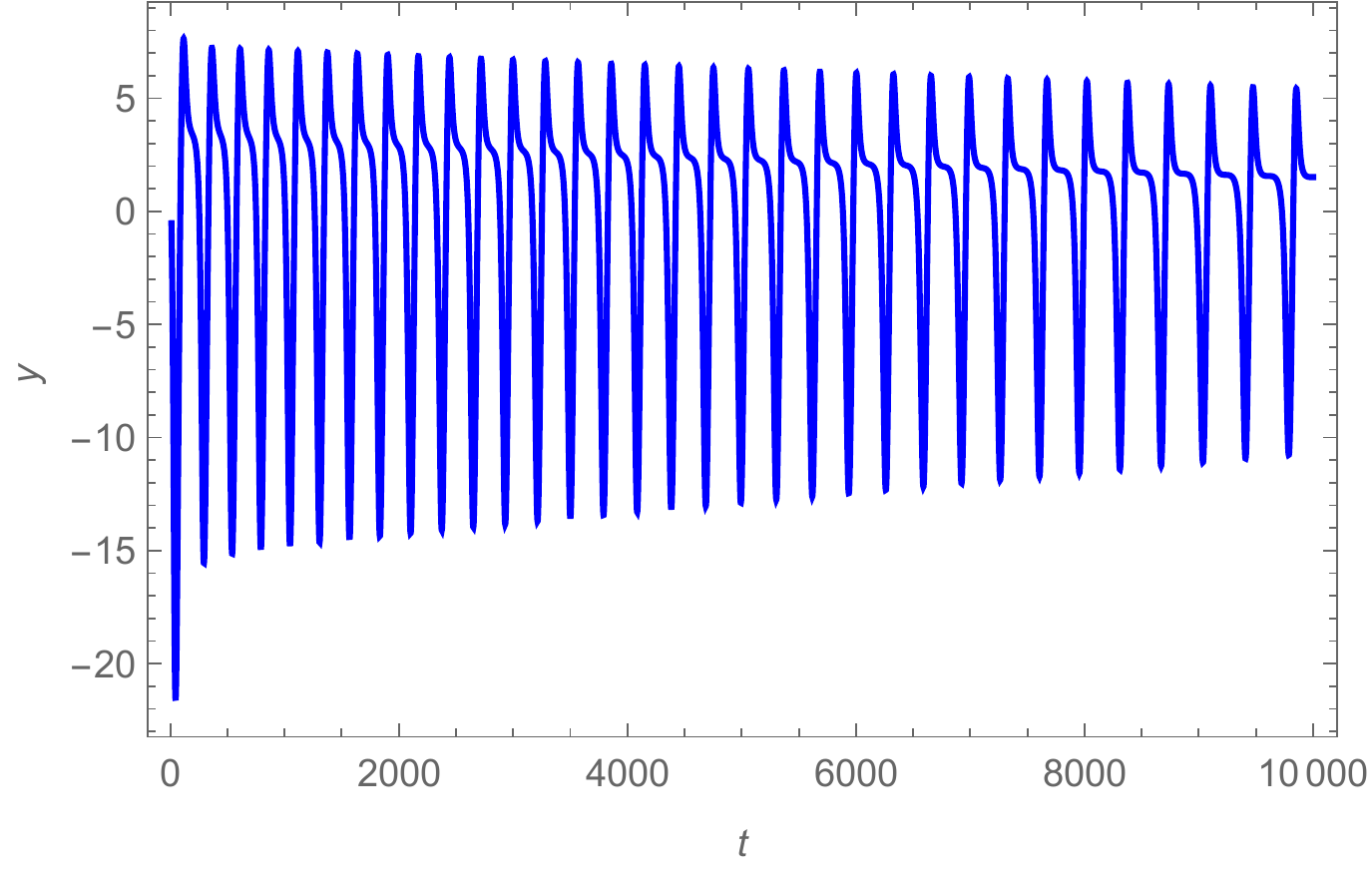}
\includegraphics[width=0.40\textwidth]{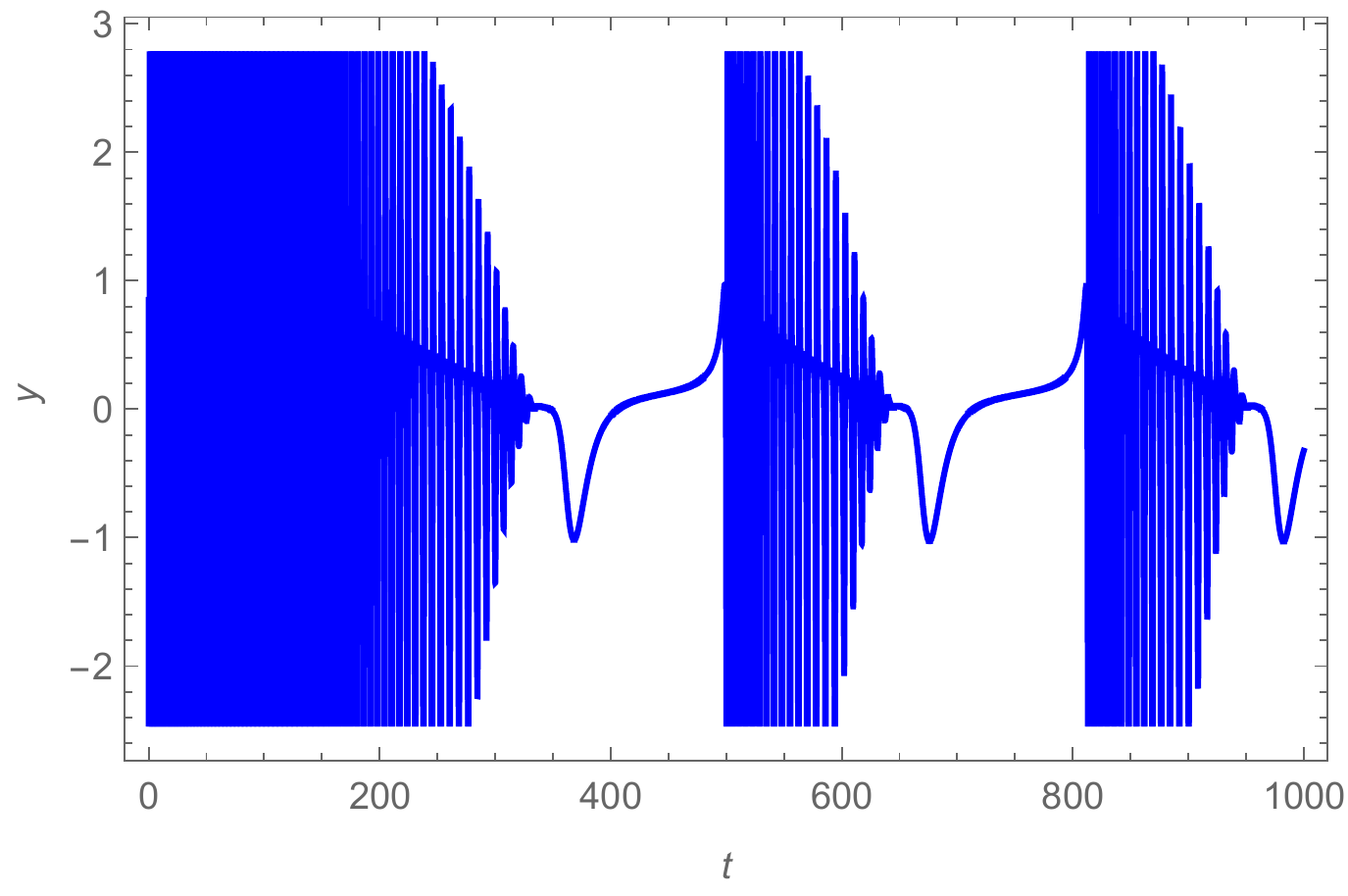}
\includegraphics[width=0.40\textwidth]{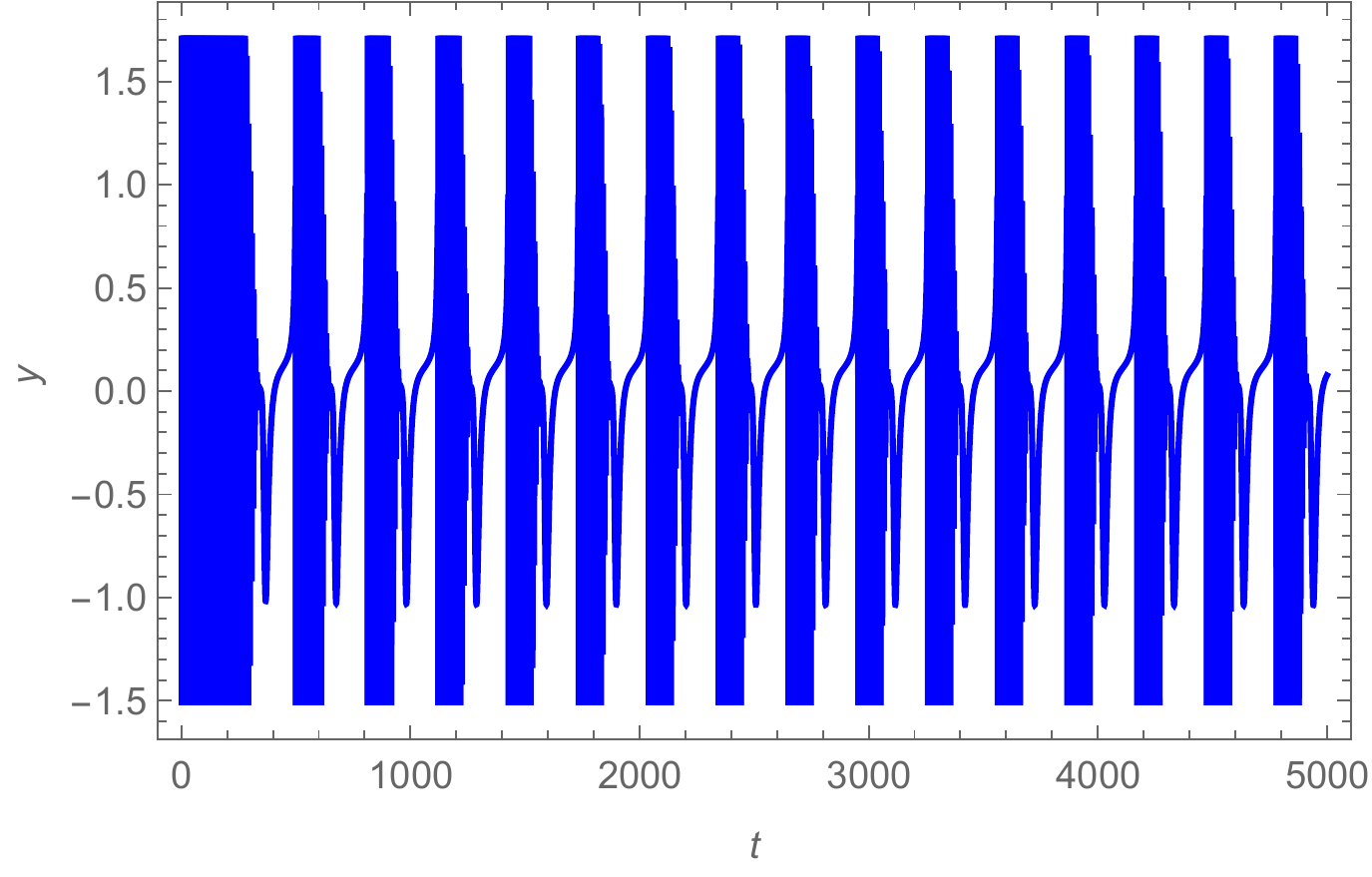}
\includegraphics[width=0.40\textwidth]{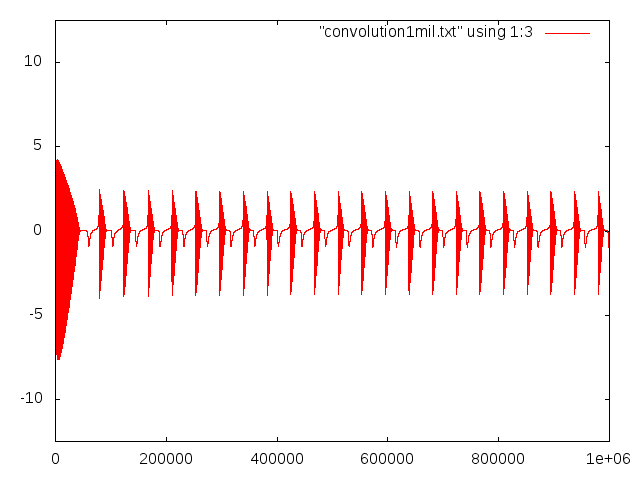}
\end{center}
\vspace{-0.5cm}
\caption{3D Hindmarsh-Rose model on 10K, 100K, 500K, 1000K iterations}
\label{fig:numerical}
\end{figure}

In both approaches we simulate the three-dimensional Hindmarsh-Rose model \cite{Hindmarsh1984} and we observed (Figure \ref{fig:numerical}) that is important to have a big number of steps in numerical simulation because the solution stability cannot be seen even after a huge number of iterations, which could take very long time to simulate without parallel processing.

\section{Conclusions and future work}
 
We observed that fixing the number of processes that we need for a run is a huge problem which can be solved by adding new processes as the iteration progress or with iterative algorithms that does not start all the treads/processes at the very early beginning.

Using Cuda to simulate such an algorithm was not a bad idea, however computing at each step of evolution the convolution for the entire history needs complicated reductions algorithms, which should be investigated closely in the future.

We present here few results from our work in progress, as a future work we should continue to investigate in both HPC and GPU approaches, what is the best way to find the numerical solution for such fractional-order systems. Further simulations would possible give us a way for the best compromise in both cases HPC od GPU.

\section{Acknowledgements}
This research has been supported by CNCS IDEI Grant PN-II-RU-TE-2014-4 - FraDys: "Theoretical and Numerical Analysis of Fractional-Order Dynamical Systems and Applications".


\begin{thebibliography}{}
{\footnotesize
\bibitem[1]{Barbosa2007}
Ramiro S. Barbosa, J.A. Tenreiro Machado, B.M. Vinagre, and A.J. Calderon.
\newblock Analysis of the Van der Pol oscillator containing derivatives of fractional order. 
\newblock {\it Journal of Vibration and Control}, 13(9-10):1291-1301, 2007.


\bibitem[2]{Donato2008Bif}
Donato Cafagna and Giuseppe Grassi.
\newblock Bifurcation and chaos in the fractional-order Chen system via a time-domain approach.
\newblock {\it International Journal of Bifurcation and Chaos}, 18(7):1845-1863, 2008.

\bibitem[3]{Donato2008Frac}
Donato Cafagna and Giuseppe Grassi.
\newblock Fractional-order Chua's circuit: time-domain analysis, bifurcation, chaotic behavior and test for chaos.
\newblock {\it International Journal of Bifurcation and Chaos}, 18(3):615-639, 2008.

\bibitem[4]{Weiwei2009}
Weiwei Zhang, Shangbo Zhou, Hua Li, and Hao Zhu.
\newblock Chaos in a fractional-order Rossler system.
\newblock {\it Chaos, Solutions and Fractals}, 42(3):1684-1691, 2009.

\bibitem[5]{Wang2010}
Mohammed Salah Abd-Elouahab, Nasr-Eddine Hamri, and Junwei Wang.
\newblock Chaos control of a fractional-order financial system.
\newblock {\it Mathematical Problems in Engineering}, 2010(Article ID 270646):1-18, 2010.

\bibitem[6]{Diethelm2002}
K. Diethelm, N.J. Ford, and A.D. Freed.
\newblock A predictor-corrector approach for the numerical solution of fractional differential equations.
\newblock {\it Nonlinear Dynamics}, 29(1-4):3-22, 2002.

\bibitem[7]{Galeone2009}
Luciano Galeone and Roberto Garrappa.
\newblock Explicit methods for fractional differential equations and their stability properties.
\newblock {\it Journal of Computational and Applied Mathematics}, 228(2):548-560, 2009.

\bibitem[8]{Diethelm2011}
Diethelm, K. 
\newblock An efficient parallel algorithm for the numerical solution of fractional differential equations.
\newblock {\it Fractional Calculus and Applied Analysis}, 14 (3): 475-490, 2011.

\bibitem[9]{Hindmarsh1984} J. Hindmarsh, R. Rose, 
\newblock A model of neuronal bursting using three coupled first order differential equations, 
\newblock {\it Proceedings of the Royal Society of London}, B221: 87-102, 1984.
}
\end{thebibliography}
\end{document}